\documentclass[12pt,reqno]{amsart}
\usepackage{amssymb}
\usepackage{amsxtra}
\usepackage[mathscr]{eucal}

\DeclareMathOperator*{\Cat}{\mathbf{Cat}}

\newcommand{\qbinom}[2]{\begin{bmatrix}#1 \\ #2\end{bmatrix}}

\newcommand{\Z}{{\mathbf Z}}

\newcommand{\Q}{{\mathbf Q}}

\newcommand{\eop}{\hfill$\square$}

\theoremstyle{plain}
\newtheorem{Thm}{Theorem}
\newtheorem{Cor}{Corollary}
\newtheorem{Lem}{Lemma}
\newtheorem{Prop}{Proposition}
\theoremstyle{definition}
\newtheorem{Def}{Definition}
\theoremstyle{remark}
\newtheorem{Rem}{Remark}

\newtheorem{Ex}{Example}

\numberwithin{equation}{section} \interdisplaylinepenalty=500

\begin{document}

\title[Multiple Harmonic $q$-Series]{Duality for Finite
Multiple Harmonic $q$-Series}

\date{\today}

\author{David~M. Bradley}
\address{Department of Mathematics \& Statistics\\
         University of Maine\\
         5752 Neville Hall
         Orono, Maine 04469-5752\\
         U.S.A.}
\email[]{bradley@math.umaine.edu, dbradley@member.ams.org}



\begin{abstract} We define two finite $q$-analogs of certain multiple
harmonic series with an arbitrary number of free parameters, and
prove identities for these $q$-analogs, expressing them in terms
of multiply nested sums involving the Gaussian binomial
coefficients.  Special cases of these identities---for example,
with all parameters equal to 1---have occurred in the literature.
The special case with only one parameter reduces to an identity
for the divisor generating function, which has received some
attention in connection with problems in sorting theory.  The
general case can be viewed as a duality result, reminiscent of the
duality relation for the ordinary multiple zeta function.
\end{abstract}

\maketitle

\section{Introduction}\label{sect:Intro}
One of the main results in~\cite{Uchi81} is the identity
\begin{equation}
   \sum_{k=1}^\infty \frac{(-1)^{k+1} q^{k(k+1)/2}}
   {(1-q)(1-q^2)\cdots(1-q^{k-1})(1-q^k)^2}
   = \sum_{k=1}^\infty \frac{q^k}{1-q^k},
   \qquad |q|<1.
\label{Keisuke}
\end{equation}
As Uchimura showed, the series~\eqref{Keisuke} is the limiting
case $(n\to\infty)$ of the polynomials
\[
   U_n(q) := \sum_{k=1}^n kq^k \prod_{j=k+1}^n (1-q^j)
\]
which arise in sorting: $U_n(1/2)$ is the average number of
exchanges required to insert a new element into a heap-ordered
binary tree with $2^n-1$ elements~\cite{Uchi77}.  Andrews and
Uchimura~\cite{AndUchi} later proved the finite analog
\begin{equation}
   \sum_{k=1}^n \frac{(-1)^{k+1} q^{k(k+1)/2}}{1-q^k}
   \qbinom{n}{k}
   = \sum_{k=1}^n \frac{q^k}{1-q^k},
   \qquad 0< n\in\Z,
\label{George}
\end{equation}
of~\eqref{Keisuke} by differentiating the $q$-Chu-Vandermonde sum.
Here,
\begin{equation}
   \qbinom{n}{k} := \prod_{j=1}^k \frac{1-q^{n-k+j}}{1-q^j}
\label{Defqbinom1}
\end{equation}
is the (Gaussian) $q$-binomial coefficient, which vanishes by
convention unless $0\le k\le n$.  Subsequently,
Dilcher~\cite{Dil95} proved the generalization
\begin{equation}
   \sum_{k=1}^n \frac{(-1)^{k+1}q^{k(k+1)/2+(m-1)k}}{(1-q^k)^m}
   \qbinom{n}{k}
   = \sum_{k_1=1}^n \frac{q^{k_1}}{1-q^{k_1}}\sum_{k_2=1}^{k_1}
   \frac{q^{k_2}}{1-q^{k_2}}\cdots\sum_{k_m=1}^{k_{m-1}}\frac{q^{k_m}}
   {1-q^{k_m}}
\label{qKarl}
\end{equation}
of~\eqref{George} by a double induction on $n$ and $m$, and noted
that~\eqref{qKarl} can be viewed as a $q$-analog of the identity
\begin{equation}
   \sum_{k=1}^n \frac{(-1)^{k+1}}{k^m}\binom{n}{k}
   =\sum_{n\ge k_1\ge k_2\ge\cdots\ge k_m\ge 1}\;\frac{1}{k_1k_2\cdots k_m},
   \qquad 0< m\in\Z,
\label{Karl}
\end{equation}
involving a finite multiple harmonic sum with unit exponents.
In this paper, we provide a generalization of~\eqref{qKarl} that
gives a $q$-analog of~\eqref{Karl} with arbitrary positive integer
exponents on the $k_j$.  We also give a $q$-analog of~\eqref{Karl}
in which the inequalities on the indices are strict.

\section{Main Result}\label{sect:Main}

Henceforth, we assume $q$ is real and $0<q<1$.  The $q$-analog of
a non-negative integer $n$ is
\[
   [n]_q := \sum_{k=0}^{n-1} q^k = \frac{1-q^n}{1-q}.
\]
\begin{Def} Let $n$, $m$ and $s_1,s_2,\dots,s_m$ be non-negative
integers.  Define the multiply nested sums
\begin{align}
   Z_n[s_1,\dots,s_m] &:= \sum_{n\ge k_1\ge \cdots\ge k_m\ge 1}
   \; \prod_{j=1}^m q^{k_j}[k_j]_q^{-s_j},
   \label{Znqdef}\\
   A_n[s_1,\dots,s_m] &:= \sum_{n\ge k_1\ge \cdots\ge k_m\ge 1}
   \; (-1)^{k_1+1} q^{k_1(k_1+1)/2}\qbinom{n}{k_1}\prod_{j=1}^m
   q^{(s_j-1)k_j} [k_j]_q^{-s_j},
   \label{Anqdef}
\end{align}
with the understanding that
$Z_0[s_1,\dots,s_m]=A_0[s_1,\dots,s_m]=0$, and with empty argument
lists, $Z_n[\;]=A_n[\;]=1$ if $n>0$ and $m=0$.  As
in~\eqref{qKarl} and~\eqref{Karl}, the sums are over all integers
$k_j$ satisfying the indicated inequalities.
\end{Def}
It will be convenient to make occasional use of the abbreviations
$\Cat_{j=1}^m \{s_j\}$ for the concatenated argument sequence
$s_1,\dots,s_m$ and $\{s\}^m=\Cat_{j=1}^m \{s\}$ for $m\ge 0$
consecutive copies of the argument $s$. We can now state our main
result.

\begin{Thm}\label{thm:main} Let $n,r$ and $a_1,b_1,\dots,a_r,b_r$ be positive
integers.  Then
\[
   Z_n[\Cat_{j=1}^{r-1}
   \{\{1\}^{a_j-1},b_j+1\},\{1\}^{a_r-1},b_r]
   = A_n[a_1,\{1\}^{b_1-1},\Cat_{j=2}^r \{a_j+1,\{1\}^{b_j-1}\}].
\]
\end{Thm}

\begin{Ex}  Putting $r=2$, $a_1=3$, $b_1=2$, $a_2=1$, $b_2=1$ in
Theorem~\ref{thm:main} gives the identity
$Z_n[1,1,3,1]=A_n[3,1,2]$, or equivalently,
\[
   \sum_{n\ge j\ge k\ge m\ge p\ge 1}
   \frac{q^{j+k+m+p}}{[j]_q[k]_q[m]_q^3[p]_q}
   =\sum_{n\ge k\ge m\ge p\ge 1}\;
   (-1)^{k+1}q^{k(k+1)/2}\qbinom{n}{k}\frac{q^{2k+p}}{[k]_q^3[m]_q[p]_q^2}.
\]
\end{Ex}

\begin{Ex} Putting $r=2$, $a_1=1$, $b_1=1$, $a_2=1$, $b_2=2$ in
Theorem~\ref{thm:main} gives the identity $Z_n[2,2]=A_n[1,2,1]$,
or equivalently,
\[
   \sum_{n\ge k\ge m\ge
   1}\frac{q^{k+m}}{[k]_q^2[m]_q^2}
   =\sum_{n\ge k\ge m\ge p\ge 1}\;
   (-1)^{k+1}q^{k(k+1)/2}\qbinom{n}{k}\frac{q^m}{[k]_q[m]_q^2[p]_q}.
\]
\end{Ex}

Theorem~\ref{thm:main} has a concise reformulation in terms of
involutions on sequences, or equivalently, dual words in the
non-commutative polynomial algebra $\Q\langle x,y\rangle$: see
Section~\ref{sect:dual}.  Additional consequences of
Theorem~\ref{thm:main} will be explored in the next section.

\section{Special Cases}
For real $x$ and $y$, we depart from convention and borrow the
suggestive notation of~\cite{QCalc} for the $q$-analog of
$(x+y)^n$:
\[
   (x+y)_q^n := \prod_{k=0}^{n-1} (x+yq^k), \qquad 0\le n\in\Z.
\]
It is easily seen that the $q$-binomial
coefficient~\eqref{Defqbinom1} has the alternative representation
\[
   \qbinom{n}{k} = \frac{(1-q)_q^n}{(1-q)_q^k (1-q)_q^{n-k}},
   \qquad 0\le k\le n,
\]
from which follow the well-known limiting results
\[
   \lim_{n\to\infty}\qbinom{n}{k} = \frac{1}{(1-q)_q^k},
   \qquad
   \lim_{q\to 1}\qbinom{n}{k} = \binom{n}{k} =
   \frac{n!}{k!(n-k)!}.
\]
Notation for limiting cases of the sums in Definition 1 are as
follows.

\begin{Def}
\begin{align*}
   Z[s_1,\dots,s_m] &:= \lim_{n\to\infty} Z_n\big[\Cat_{j=1}^m s_j\big]
   =\sum_{k_1\ge \cdots\ge k_m\ge 1}
   \; \prod_{j=1}^m q^{k_j}[k_j]_q^{-s_j},
   \\
   A[s_1,\dots,s_m] &:= \lim_{n\to\infty} A_n\big[\Cat_{j=1}^m s_j\big]
   =\sum_{k_1\ge \cdots\ge k_m\ge 1}
   \; \frac{(-1)^{k_1+1} q^{k_1(k_1+1)/2}}{(1-q)_q^{k_1}}\prod_{j=1}^m
   \frac{q^{(s_j-1)k_j}}{[k_j]_q^{s_j}}.
\end{align*}
\end{Def}

\begin{Def}\label{ZAdefs}
\begin{align*}
   Z_n(s_1,\dots,s_m) &:= \lim_{q\to 1} Z_n\big[\Cat_{j=1}^m s_j\big]
   =\sum_{n\ge k_1\ge \cdots\ge k_m\ge 1}
   \; \prod_{j=1}^m k_j^{-s_j},
   \\
   A_n(s_1,\dots,s_m) &:= \lim_{q\to 1} A_n\big[\Cat_{j=1}^m s_j\big]
   = \sum_{n\ge k_1\ge\cdots\ge k_m\ge 1}\;
    (-1)^{k_1+1}\binom{n}{k_1}\prod_{j=1}^m k_j^{-s_j}.
\end{align*}
\end{Def}

With this notation, the following consequences of Theorem 1 are
immediate.

\begin{Cor}\label{cor:limiting} Let $n,r$ and $a_1,b_1,\dots,a_r,b_r$ be positive
integers.  Then
\begin{align*}
   Z_n(\Cat_{j=1}^{r-1}
   \{\{1\}^{a_j-1},b_j+1\},\{1\}^{a_r-1},b_r)
   &= A_n(a_1,\{1\}^{b_1-1},\Cat_{j=2}^r \{a_j+1,\{1\}^{b_j-1}\}),
   \\
   Z[\Cat_{j=1}^{r-1}
   \{\{1\}^{a_j-1},b_j+1\},\{1\}^{a_r-1},b_r]
   &= A[a_1,\{1\}^{b_1-1},\Cat_{j=2}^r \{a_j+1,\{1\}^{b_j-1}\}].
\end{align*}
\end{Cor}

\begin{Cor}\label{cor:qKarl} Let $n$, $a$ and $b$ be positive integers.
Then
\begin{equation}
   Z_n[\{1\}^{a-1},b] = A_n[a,\{1\}^{b-1}].
\label{duality}
\end{equation}
\end{Cor}
\noindent{\bf Proof.} Put $r=1$ in Theorem~\ref{thm:main}. \eop

\begin{Rem} If we put $b=1$ and $a=m$ in~\eqref{duality}, we find that
\begin{equation}
   Z_n[\{1\}^m]=A_n[m],
\label{An[m]}
\end{equation}
which is~\eqref{qKarl}.  As we shall see, $A_n[0]=1$ is an easy
consequence of the $q$-binomial theorem.  Since we have defined
$Z_n[\;]=1$, it follows that~\eqref{An[m]} also holds when $m=0$.
Concerning the equivalent equation~\eqref{qKarl}, this point was
also noted by Dilcher~\cite{Dil95}.
\end{Rem}

\begin{Rem}  If we put $a=1$ and $b=m$
in~\eqref{duality}, there follows $Z_n[m]=A_n[\{1\}^m]$, an
identity dual to~\eqref{qKarl} and~\eqref{An[m]}:
\begin{equation}\label{qKarldual}
   \sum_{k=1}^n \frac{q^k}{[k]_q^m} =
   \sum_{n\ge k_1\ge\cdots\ge k_m\ge 1}\;
   (-1)^{k_1+1} q^{k_1(k_1+1)/2}\qbinom{n}{k_1}\prod_{j=1}^m
   \frac{1}{[k_j]_q},
\end{equation}
with respective limiting cases
\[
   \sum_{k=1}^n \frac{1}{k^m} =
   \sum_{n\ge k_1\ge\cdots\ge k_m\ge 1}\;
   (-1)^{k_1+1}\binom{n}{k_1}\prod_{j=1}^m
   \frac{1}{k_j}
\]
and
\[
   \sum_{k=1}^\infty \frac{q^k}{[k]_q^m}
   = \sum_{k_1\ge \cdots\ge k_m\ge 1}\;
   \frac{(-1)^{k_1+1}q^{k_1(k_1+1)/2}}{(1-q)_q^{k_1}}
   \prod_{j=1}^m \frac{1}{[k_j]_q}.
\]
Prodinger's main result~\cite[Theorem 1]{Prod} is easily obtained
from~\eqref{qKarldual} by replacing $q$ with $1/q$. The fact that
Prodinger obtained his result from~\eqref{qKarl} as a consequence
of an inverse-pair equivalence suggests that additional instances
of our Theorem~\ref{thm:main} may likewise be so related.  That
this is indeed the case is one of the insights of the next
section.
\end{Rem}

\section{Duality}\label{sect:dual}
Here, we recast Theorem~\ref{thm:main} in the language of
\emph{duality}, a concept first formulated for multiple harmonic
series in~\cite{Hoff1}, and subsequently
generalized~\cite{Ohno,DBqMzv,Okuda}.  Following a suggestion of
Hoffman~\cite{Hoff0}, define an involution on the set $\mathscr S$
of finite sequences of positive integers as follows. Let $\alpha$
be the map that sends a sequence in $\mathscr{S}$ to its sequence
of partial sums. The image of $\alpha$ thus consists of the
strictly increasing finite sequences of positive integers, on
which we define an involution $\beta$ by
\[
   \beta\big(\Cat_{j=1}^m t_j\big) =
   \{k\in\Z : 1\le k\le t_m \}\setminus\{t_j : 1\le j\le m-1\},
   \qquad 0<m\in\Z,
\]
arranged in increasing order.  In other words, $\beta$ maps a
strictly increasing sequence of positive integers $t_1,\dots,t_m$
to its set-theoretic complement in the positive integers up to
$t_m$, and then tacks $t_m$ onto the end of the result. Clearly,
the composition of maps $\alpha^{-1}\beta\alpha$ is an involution
of $\mathscr S$, and it is easy to see that Theorem~\ref{thm:main}
can be restated as
\[
   Z_n[\vec s\,] = A_n[\alpha^{-1}\beta\alpha \vec s\,],
   \qquad \forall \vec s\in\mathscr{S},\quad 0<n\in\Z.
\]

For an alternative duality reformulation, let $\mathfrak h =
\Q\langle x,y\rangle$ denote the non-commutative polynomial
algebra over the field of rational numbers in two indeterminates
$x$ and $y$.
Let $\mathfrak h' = \mathfrak hy$ and fix a positive integer $n$.
Define $\Q$-linear maps $\widehat{A}_n$ and $\widehat{Z}_n$ on
$\mathfrak h'$ by
\[
   \widehat{A}_n\bigg[\prod_{j=1}^m x^{s_j-1}y\bigg]
   := A_n\big[\Cat_{j=1}^m s_j\big], \qquad
   \widehat{Z}_n\bigg[\prod_{j=1}^m x^{s_j-1}y\bigg]
   := Z_n\big[\Cat_{j=1}^m s_j\big],
\]
for any positive integers $s_1,s_2,\dots,s_m$. Let $J$ be the
automorphism of $\mathfrak h$ that switches $x$ and $y$. Define an
involution of $\mathfrak h'$ by
\begin{equation}\label{dualword}
   w^* = (Jw)x^{-1}y,\qquad \forall w\in \mathfrak h'.
\end{equation}
It is now a routine matter to show that Theorem~\ref{thm:main} can
be restated as
\begin{equation}\label{AZdual}
  \widehat{A}_n[w] = \widehat{Z}_n[w^*],\qquad
  \forall w\in\mathfrak h'.
\end{equation}
Now Prodinger~\cite[Lemma 1]{Prod} proved that for positive
integer $n$, the inverse pairs
\begin{align}
   \sum_{k=0}^n \beta_k
    &= \sum_{k=0}^n (-1)^k q^{k(k-1)/2}\qbinom{n}{k} \alpha_k\label{pro1}\\
\intertext{and}
   \sum_{k=0}^n q^{-k} \alpha_k
    &= \sum_{k=0}^n (-1)^k q^{k(k-1)/2-kn}\qbinom{n}{k} \beta_k\label{pro2}
\end{align}
are equivalent.  In other words,~\eqref{pro1} holds for a pair of
sequences $\alpha_0,\alpha_1,\dots,\alpha_n$ and
$\beta_0,\beta_1,\dots,\beta_n$ if and only if~\eqref{pro2} does.
But our Theorem~\ref{thm:main} states that~\eqref{pro1} holds for
the sequences
\[
   \alpha_k = -\sum_{k\ge k_2\ge\cdots\ge k_m\ge 1}\;
   \frac{q^{s_1 k}}{[k]_q^{s_1}} \prod_{j=2}^m
   \frac{q^{(s_j-1)k_j}}{[k_j]_q^{s_j}},\qquad
   \beta_k = \sum_{k\ge k_2\ge\cdots\ge k_m\ge 1}\;
   \frac{q^k}{[k]_q^{s_1}}\prod_{j=2}^m
   \frac{q^{k_j}}{[k_j]_q^{s_j}},
\]
where $s_1,\dots,s_m$ are positive integers, $1\le k\le n$, and
$\alpha_0=\beta_0=0$.  Equivalently,~\eqref{AZdual} holds with
$w=x^{s_1-1}y\cdots x^{s_m-1}y$.  Since duality is an involution
on words, $w^{**}=w$, and Theorem~\ref{thm:main} also gives the
dual statement $A_n[w^*]=Z_n[w]$. But this latter statement is
easily seen to be equivalent to~\eqref{pro2} if we replace $q$ by
$1/q$ throughout. Thus Prodinger's result implies that if
$A_n[w]=Z_n[w^*]$ is known for a particular word $w$, then we also
know $A_n[w^*]=Z_n[w]$, and vice-versa.  Of course, knowing only
that the two statements are equivalent does not establish their
truth in any particular instances---for this we need our
Theorem~\ref{thm:main}.  What it shows is that the notion of
duality in the case of Theorem~\ref{thm:main} coincides with the
existence of a certain class of inverse pair identities.

For additional duality results concerning multiple harmonic
$q$-series, see~\cite{DBqMzv,Okuda}.  It is interesting to
contrast~\eqref{dualword} and~\eqref{AZdual} with the
corresponding duality statement for multiple zeta
values~\cite{HoffAlg,HoffOhno}, for which the relevant involution
is the \emph{anti}-automorphism of $x\mathfrak hy$ that switches
$x$ and $y$.

\section{Proof of Theorem~\ref{thm:main}}\label{sect:proof}
By induction, it suffices to establish the two recurrence
relations for $A_n$ stated in Propositions~\ref{prop:AnSum}
and~\ref{prop:An0} below, along with the base cases
$A_n[\;]=A_n[0]=1$ for $0<n\in\Z$.

\begin{Prop}\label{prop:AnSum} Let $n$, $m$, and $s_1,s_2,\dots,s_m$ be positive
integers. Then
\[
   A_n[s_1,s_2,\dots,s_m] = \sum_{r=1}^n q^r [r]_q^{-1}
   A_r[s_1-1,s_2,\dots,s_m].
\]
\end{Prop}

\begin{Prop}\label{prop:An0} Let $n$, $m$, and $s_2,s_3,\dots,s_m$ be
positive integers.  Then
\[
   A_n[0,s_2,s_3,\dots,s_m] = [n]_q^{-1} A_n[s_2-1,s_3,\dots,s_m].
\]
\end{Prop}

The base case $A_n[\;]=1$ for $n>0$ is true by definition. As
alluded to previously, the other base case is an easy consequence
of the $q$-binomial theorem~\cite{AndAskRoy,Gasp,QCalc}
\begin{equation}
   (x+y)_q^n = \sum_{m=0}^n q^{m(m-1)/2}\qbinom{n}{m}
   x^{n-m} y^m.
\label{qbinom}
\end{equation}
Putting $x=1$ and $y=-1$ in~\eqref{qbinom}, we see that if $n>0$,
then
\[
   A_n[0]
   = \sum_{m=1}^n (-1)^{m+1} q^{m(m-1)/2}\qbinom{n}{m}
   = 1-(1-1)_q^n
   = 1.
\]
Thus, it remains only to prove Propositions~\ref{prop:AnSum}
and~\ref{prop:An0}.

We shall make use of the (equivalent by symmetry) $q$-Pascal
recurrences~\cite{QCalc}
\begin{equation}
   \qbinom{r}{k} = \qbinom{r-1}{k} +
   q^{r-k}\qbinom{r-1}{k-1},\qquad\quad
   \qbinom{r}{k} =
   q^k\qbinom{r-1}{k}+\qbinom{r-1}{k-1}
   \label{qPascal}
\end{equation}
and the following elementary summation formula.

\begin{Lem}\label{lem:qsum} Let $k$ and $n$ be positive integers with $k\le n$.  Then
\[
   \sum_{r=k}^n q^r\qbinom{r-1}{k-1} = q^k\qbinom{n}{k}.
\]
\end{Lem}
\noindent{\bf Proof.}  Write the first $q$-Pascal
recurrence~\eqref{qPascal} in the form
\[
   q^{r-k}\qbinom{r-1}{k-1} = \qbinom{r}{k}-\qbinom{r-1}{k},
\]
multiply through by $q^k$, and sum on $r$. \eop

It will be helpful to introduce some further notation.
\begin{Def}\label{Wndef} Let $n$, $m$, and $s_1,\dots,s_m$ be non-negative
integers.  Define
\[
   W_n[s_1,\dots,s_m] := \sum_{n\ge k_1\ge\cdots\ge k_m\ge 1}
   \; \prod_{j=1}^m q^{(s_j-1)k_j}[k_j]_q^{-s_j},
\]
with the understanding that $W_0[s_1,\dots,s_m]=0$ and $W_n[\;]=1$
if $n>0$ and $m=0$.
\end{Def}

\noindent{\bf Proof of Proposition~\ref{prop:AnSum}.}  By
Lemma~\ref{lem:qsum}, we have
\begin{align*}
   A_n[s_1,\dots,s_m] &= \sum_{k=1}^n (-1)^{k+1}
   q^{k(k-1)/2+(s_1-1)k}[k]_q^{-s_1}W_k[s_2,\dots,s_m]
   q^k\qbinom{n}{k}\\
   &=\sum_{k=1}^n (-1)^{k+1}
   q^{k(k-1)/2+(s_1-1)k}[k]_q^{-s_1}W_k[s_2,\dots,s_m]
   \sum_{r=k}^n q^r\qbinom{r-1}{k-1}\\
   &=\sum_{r=1}^n q^r\sum_{k=1}^r
   (-1)^{k+1}q^{k(k-1)/2+(s_1-1)k}\qbinom{r-1}{k-1}
   [k]_q^{-s_1}W_k[s_2,\dots,s_m]\\
   &=\sum_{r=1}^n q^r[r]_q^{-1}\sum_{k=1}^r (-1)^{k+1}
   q^{k(k+1)/2+(s_1-2)k}\qbinom{r}{k}[k]_q^{1-s_1}W_k[s_2,\dots,s_m]\\
   &=\sum_{r=1}^n q^r[r]_q^{-1}A_r[s_1-1,s_2,\dots,s_m].
\end{align*}
\eop

\noindent{\bf Proof of Proposition~\ref{prop:An0}.}  By the second
$q$-Pascal recurrence~\eqref{qPascal}, we have
\begin{align*}
   A_n[0,s_2,\dots,s_m] &=\sum_{k=1}^n (-1)^{k+1}
   q^{k(k-1)/2}\qbinom{n}{k}W_k[s_2,\dots,s_m]\\
   &= \sum_{k=1}^n (-1)^{k+1} q^{k(k-1)/2}
   \bigg\{q^k\qbinom{n-1}{k}+\qbinom{n-1}{k-1}\bigg\}W_k[s_2,\dots,s_m]\\
   &=\sum_{k=1}^n
   (-1)^{k+1}q^{k(k+1)/2}\qbinom{n-1}{k}W_k
   +\sum_{k=1}^n (-1)^{k+1} q^{k(k-1)/2}\qbinom{n-1}{k-1}W_k\\
   &=\sum_{k=1}^{n-1}
   (-1)^{k+1}q^{k(k+1)/2}\qbinom{n-1}{k}W_k
   +\sum_{k=0}^{n-1} (-1)^k q^{k(k+1)/2}\qbinom{n-1}{k}W_{k+1}\\
   &=W_1+\sum_{k=1}^{n-1} (-1)^k
   q^{k(k+1)/2}\qbinom{n-1}{k}(W_{k+1}-W_k)\\
   &= \sum_{k=0}^{n-1}(-1)^k
   q^{k(k+1)/2}\qbinom{n-1}{k}q^{(s_2-1)(k+1)}\,[k+1]_q^{-s_2}\,W_{k+1}[s_3,\dots,s_m]\\
   &=\sum_{k=1}^n (-1)^{k+1}
   q^{k(k-1)/2}\qbinom{n-1}{k-1}q^{(s_2-1)k}\,[k]_q^{-s_2}\,W_k[s_3,\dots,s_m]\\
   &=[n]_q^{-1}\sum_{k=1}^n (-1)^{k+1}
   q^{k(k+1)/2}\qbinom{n}{k}q^{(s_2-2)k}\,[k]_q^{1-s_2}\,W_k[s_3,\dots,s_m]\\
   &= [n]_q^{-1}A_n[s_2-1,s_3,\dots,s_m].
\end{align*}
\eop

The following consequence of Proposition~\ref{prop:An0} may be
worth noting.
\begin{Cor}\label{cor:An01m} Let $n$ be a positive integer, and let $m$ be a non-negative integer.
Then $A_n[0,\{1\}^m] = [n]_q^{-m}$, or equivalently,
\[
   \sum_{n\ge k\ge k_1\ge\cdots\ge k_m\ge 1} (-1)^{k+1} q^{k(k-1)/2} \qbinom{n}{k}
   \prod_{j=1}^m [k_j]_q^{-1} = [n]_q^{-m}.
\]
\end{Cor}

\section{Multiple Harmonic $q$-Series With Strict Inequalities}
Much of the recent literature concerning multiple harmonic series
has focused on sums of the form
\begin{equation}
   \zeta(s_1,\dots,s_m) := \sum_{k_1>\cdots>k_m>0}\;\prod_{j=1}^m
   k_j^{-s_j}
\label{mzvdef}
\end{equation}
and various multiple polylogarithmic
extensions~\cite{BBBLa,BBBLc,BowBrad1,BowBrad2,BowBrad3,BowBradRyoo,prtn,Hoff1,HoffAlg,HoffOhno,Ohno}.
Thus, it may be of interest to consider finite $q$-analogs
of~\eqref{mzvdef} akin to~\eqref{Znqdef}, but in which the
inequalities are strict as opposed to weak.
\begin{Def}\label{def:strict}  Let $n$, $m$ and $s_1,s_2,\dots,s_m$ be non-negative
integers.  Define
\begin{align*}
   Z_n^{>}[s_1,\dots,s_m] &:= \sum_{n\ge k_1> \cdots > k_m\ge 1}
   \; \prod_{j=1}^m q^{k_j}[k_j]_q^{-s_j},
   \\
   A_n^{>}[s_1,\dots,s_m] &:= \sum_{n\ge k_1>\cdots > k_m\ge 1}
   \; (-1)^{k_1} q^{k_1(k_1+1)/2}\qbinom{n}{k_1}\prod_{j=1}^m
   q^{(s_j-1)k_j} [k_j]_q^{-s_j},
\end{align*}
with the understanding that if $m>0$, then
$Z_0^{>}[s_1,\dots,s_m]=A_0^{>}[s_1,\dots,s_m]=0$ and if $m=0$,
then $Z_n^{>}[\;]=A_n^{>}[\;]=1$ for all $n\ge0$.
\end{Def}
In light of~\eqref{mzvdef}, it is clear that
 $  \lim_{n\to\infty}\lim_{q\to1}Z_n^{>}[\vec s\,] = \zeta(\vec s\,)$
if $\vec s=s_1,\dots,s_m$ is any vector of positive integers with
$s_1>1$.  Of course, there is an obvious relationship between
$Z_n$ and $Z_n^{>}$ and between $A_n$ and $A_n^{>}$. For example,
\[
   Z_n[s]=Z_n^{>}[s],\qquad
   Z_n[s_1,s_2]=Z_n^{>}[s_1,s_2]+Z_n^{>}[s_1+s_2],
\]
and
\[
   Z_n[s_1,s_2,s_3]=Z_n^{>}[s_1,s_2,s_3]+Z_n^{>}[s_1+s_2,s_3]+Z_n^{>}[s_1,s_2+s_3]+Z_n^{>}[s_1+s_2+s_3].
\]
More generally, $Z_n[\vec s\,]$ is the sum of those $Z_n^{>}[\vec
t\,]$, where $\vec t$ is obtained from $\vec s$ by replacing any
number of commas by plus signs.
Despite this relationship, the presence of strict inequalities in
Definition~\ref{def:strict} does appear to complicate matters
insofar as there is no simple analog of Theorem~\ref{thm:main}
relating $Z_n^{>}$ to $A_n^{>}.$ Nevertheless, there are
recurrences for $A_n^{>}$ analogous to the recurrences satisfied
by $A_n$. Arguing as in Section~\ref{sect:proof}, we find that
\begin{align*}
   A_n^{>}[s_1,\dots,s_m] &= \sum_{r=1}^n
   q^r[r]_q^{-1}A_r^{>}[s_1-1,s_2,\dots,s_m],\\
   A_n^{>}[0,s_2,\dots,s_m] &= -A_{n-1}^{>}[s_2,\dots,s_m].
\end{align*}
As a consequence, one can establish (using induction as in the
proof of Corollary~\ref{cor:An01m}) the following result.

\begin{Thm}\label{thm:Zn>ones} Let $m$ and $n$ be non-negative
integers.  Then $(-1)^mZ_n^{>}[\{1\}^m] = A_n^{>}[\{1\}^m].$
\end{Thm}

\begin{Cor} Let $m$ be a positive integer.  Then
\[
   (-1)^m\sum_{k_1>\cdots>k_m>0}\;\prod_{j=1}^m
   \frac{q^{k_j}}{[k_j]_q}=
   \sum_{k_1>\cdots>k_m>0}\;
   (-1)^{k_1}\frac{q^{k_1(k_1+1)/2}}{(1-q)_q^{k_1}}
   \prod_{j=1}^m \frac{1}{[k_j]_q}.
\]
\end{Cor}
\noindent{\bf Proof.} Let $n\to\infty$ in
Theorem~\ref{thm:Zn>ones}. \eop

\begin{Cor} Let $m$ and $n$ be positive integers.  Then
\[
   (-1)^m\sum_{n\ge k_1>\cdots >k_m\ge 1} \;
   \prod_{j=1}^m k_j^{-1} =
   \sum_{n\ge k_1>\cdots >k_m\ge 1}\;
   (-1)^{k_1}\binom{n}{k_1}\prod_{j=1}^m k_j^{-1}.
\]
\end{Cor}
\noindent{\bf Proof.} Let $q\to 1$ in
Theorem~\ref{thm:Zn>ones}.\eop

One can also prove Theorem~\ref{thm:Zn>ones} by differentiating
the $q$-Chu-Vandermonde summation as in the
Andrews-Uchimura~\cite{AndUchi} proof of~\eqref{George}:

\noindent{\bf Alternative Proof of Theorem~\ref{thm:Zn>ones}.} Let
$x$ be real, $x\ne 1$. Write the $q$-Chu-Vandermonde
sum~\cite{Gasp,Slater} in the form
\[
   \frac{(1-xq)_q^n}{(1-q)_q^n} = 1+\sum_{k=1}^n (-1)^k
   q^{k(k+1)/2}\qbinom{n}{k}\frac{(x-1)_q^k}{(1-q)_q^k}.
\]
and differentiate both sides $m>0$ times with respect to $x$,
obtaining
\begin{multline}\label{diffqChuVdm}
   (-1)^m\frac{(1-xq)_q^n}{(1-q)_q^n}\sum_{n\ge k_1>\cdots >
   k_m\ge 1}\;\prod_{j=1}^m \frac{q^{k_j}}{1-xq^{k_j}}\\
   =\sum_{k=1}^n (-1)^kq^{k(k+1)/2}\qbinom{n}{k}
   \frac{(x-1)_q^k}{(1-q)_q^k}\sum_{k>k_1>\cdots>k_m\ge 0}\;
   \prod_{j=1}^m \frac{1}{x-q^{k_j}}.
\end{multline}
Now let $x\to1$ and note that the sum on the right hand side
of~\eqref{diffqChuVdm} vanishes if $k_m>0$.  Thus,
\[
   (-1)^m\sum_{n\ge k_1>\cdots>k_m\ge
   1}\;\prod_{j=1}^m\frac{q^{k_j}}{1-q^{k_j}} =
   \sum_{k=1}^n(-1)^k\qbinom{n}{k}\frac{q^{k(k+1)/2}}{1-q^k}
   \sum_{k>k_1\cdots>k_{m-1}>0}\prod_{j=1}^{m-1}\frac{1}{1-q^{k_j}}.
\]
Finally, multiply both sides by $(1-q)^m$ to complete the
proof.\eop

\section{Final Remarks}
In~\cite{DBqMzv,Okuda,Zhao}, the multiple $q$-zeta function
\[
   \zeta[\vec s;q] = \zeta[s_1,\dots,s_m] := \sum_{k_1>\cdots >k_m>0}\;
   \prod_{j=1}^m \frac{q^{(s_j-1)k_j}}{[k_j]_q^{s_j}}
\]
is the central object of study.  With the abbreviation
\[
   |\vec s| := \sum_{j=1}^m s_j,
\]
we note the relationship
\[
   \zeta[\vec s;q] = q^{|\vec s|}Z_{\infty}^{>}[\vec s;1/q],
\]
where $Z_{\infty}^{>}[\vec s;1/q]$ denotes the limit as
$n\to\infty$ of the $Z_n^{>}$-function of
Definition~\ref{def:strict} with $q$ replaced by $1/q$.

Fu and Lascoux~\cite{FuLas} have generalized~\eqref{qKarl} in a
different direction.  They proved that if $n$ and $m$ are positive
integers, then
\begin{multline}\label{fulas}
   \sum_{n\ge k_1\ge \cdots\ge k_m\ge 1}\;\prod_{j=1}^m
   \frac{a-bq^{k_j}}{c-zq^{k_j}}\\
    =
   \frac{c^n(1-zq/c)_q^n}{(1-q)_q^n(az-bc)^{n-1}}
   \sum_{k=1}^n \qbinom{n}{k}
   \frac{(-1)^{k-1}q^{k(k+1)/2-nk}(1-q^k)(a-bq^k)^{m+n-1}}{(c-zq^k)^{m+1}}.
\end{multline}
Letting $c=z=1$, $a=0$, $b=-1$ in~\eqref{fulas}
gives~\eqref{qKarl}.

\end{document}